\documentclass[12pt,reqno]{amsart}

\usepackage{amscd}
\usepackage{amssymb,latexsym}
\usepackage{epsfig}
\usepackage{float}
\usepackage{url}
\usepackage{graphicx,xcolor}

\newcommand{\N}{{\mathbb N}}

\newcommand{\R}{{\mathbb R}}
\renewcommand{\P}{{\mathbb P}}

\newcommand{\AAA}{{\mathcal A}}
\newcommand{\BB}{{\mathcal B}}
\newcommand{\CC}{{\mathcal C}}
\newcommand{\DD}{{\mathcal D}}

\newcommand{\FF}{{\mathcal F}}

\newcommand{\NN}{{\mathcal N}}

\newcommand{\PP}{{\mathcal P}}

\newcommand{\UU}{{\mathcal U}}

\newcommand{\www}{\widetilde}

\newcommand{\oooo}{\overline}
\newcommand{\uuuu}{\underline}

\newcommand{\pr}{{\rm pr}}

\newcommand{\supp}{{\rm supp}}
\newcommand{\paa}{\partial}

\DeclareMathOperator{\codim}{codim}

\DeclareMathOperator{\Gr}{Gr}

\DeclareMathOperator{\id}{id}

\DeclareMathOperator{\rank}{rank}

\makeatletter 
\@namedef{subjclassname@2020}{\textup{2020} Mathematics Subject Classification}
\makeatother

\begin{document}

\theoremstyle{plain}
\newtheorem{lemma}{Lemma}[section]
\newtheorem{definition/lemma}[lemma]{Definition/Lemma}
\newtheorem{theorem}[lemma]{Theorem}
\newtheorem{proposition}[lemma]{Proposition}
\newtheorem{korollar}[lemma]{Korollar}
\newtheorem{behauptungen}[lemma]{Behauptungen}
\newtheorem{conjecture}[lemma]{Conjecture}
\newtheorem{conjectures}[lemma]{Conjectures}
\newtheorem{corollary}[lemma]{Corollary}

\theoremstyle{definition}
\newtheorem{definition}[lemma]{Definition}
\newtheorem{withouttitle}[lemma]{}
\newtheorem{remark}[lemma]{Remark}
\newtheorem{remarks}[lemma]{Remarks}
\newtheorem{example}[lemma]{Example}
\newtheorem{notations}[lemma]{Notations}

\title[Mixed extensions of generic finite games]
{Mixed extensions of generic finite games embedded
into products of real projective spaces} 

\author{Claus Hertling and Matija Vuji\'c}

\address{Claus Hertling\\
Lehrstuhl f\"ur algebraische Geometrie, Universit\"at Mannheim,
B6 26, 68159 Mannheim, Germany}

\email{claus.hertling@uni-mannheim.de}

\address{Matija Vuji\'c\\
Z\"urich, Switzerland}


\date{December 22, 2024}

\subjclass[2020]{91A06, 91A10}

\keywords{Generic finite games, mixed extensions,
transversality, projective spaces}

\begin{abstract}
{\Small Finite games in normal form and their mixed extensions
are a corner stone of noncooperative game theory.
Often {\it generic} finite games and their mixed extensions
are considered. But the properties which one expects in
generic games and the existence of games with these properties
are often treated only in passing. 
The paper considers strong properties and proves that generic
games have these properties.
The space of mixed strategy combinations is embedded 
in a natural way into a product of real projective spaces.
All relevant hypersurfaces extend to this bigger space.
The paper shows that for all games in the complement of a 
semialgebraic subset of codimension at least one
all relevant hypersurfaces in the bigger space are smooth and
maximally transversal. The proof uses the theorem of Sard
and follows an argument of Khovanskii.}
\end{abstract}

\maketitle

\tableofcontents

\section{Introduction}\label{c1}
\setcounter{equation}{0}

Finite games in normal form and their mixed extensions 
are a corner stone of noncooperative game theory. 
Nash showed that the mixed extension of any finite game
in normal form has Nash equilibria \cite{Na51}.
For special games, the set of Nash equilibria can be
complicated. But for generic games it is finite and odd
\cite{Wi71}\cite{Ro71}\cite{Ha73}\cite{Ri94}.
These references and many others consider {\it generic}
finite games and their mixed extensions.

The notion of a {\it generic} finite game is usually
defined ad hoc and results on it are usually only proved
as an aside. One fixes the finite set of players and for each
player his finite set of pure strategies. Then one obtains
an affine linear space $\UU$ of tuples of possible utility
functions. Generic games are then games with generic
utility functions, and these are contained in the complement
of a closed subset of measure zero.
The properties which generic games are supposed to have,
depend on the paper and the author. 

In this paper we consider a strong version of genericity
and prove that all games in the complement $\UU-\DD$
of a semialgebraic subset $\DD$ of codimension at least one
are generic. 

We fix the number $m$ of players $i\in\AAA=\{1,...,m\}$
and for each player $i$ the finite set
$S^i=\{s^i_0,...,s^i_{n_i}\}$ of $n_i+1$ pure strategies.
$S=S^1\times ...\times S^m$ is the set of
pure strategy combinations.
Then the finite game is determined by the tuple
$U=(U^1,...,U^m):S\to \R^m$ of the utility functions
$U^i:S\to\R$ of the players $i$. It is an element of the real
vector space $\UU=(\R^S)^m$.
In the mixed extension $(\AAA,G,V)$ of the finite game
$(\AAA,S,U)$, the set $G$ is
$G=G^1\times ...\times G^m$, and the set $G^i$
of mixed strategies of player $i$ is the set of
probability distributions over $S^i$, so it is a
simplex of dimension $n_i$. The utility function
$V^i:G\to\R$ is the multilinear extension of $U^i$
(see formula \eqref{2.1}).

The multilinearity leads to a natural embedding
$G^i\hookrightarrow \P^{n_i}\R$ of $G^i$ into
the real projective space $\P^{n_i}\R$ of dimension
$n_i$ and to an embedding of $G$ into the 
product $\P^\AAA W:=\prod_{i\in\AAA}\P^{n_i}\R$
of real projective spaces. 
All the induced functions, which are relevant
for understanding the game, its best reply maps
and its Nash equilibria, extend from $G$ to
the space $\P^\AAA W$. 

Our main result Theorem \ref{t3.1} says that
there is a semialgebraic subset $\DD\subset\UU$
of codimension at least one such that for each
$U\in\UU-\DD$, all natural hypersurfaces in $G$
are smooth and 
extend to smooth hypersurfaces in $\P^\AAA W$
and that they are (in a precise sense) everywhere
maximally transversal. 

This implies immediately that for such a game
the set of all Nash equilibria is finite and that
each Nash equilibrium is regular in a strong sense.

It also provides a foundation for the argument
of Wilson \cite{Wi71} and Rosenm\"uller \cite{Ro71}
that the number of Nash equilibria of a generic
game is odd. Both papers claim implicitly the
existence of generic games with good properties
without proof. The paper here provides such a proof.

Our proof follows an argument of Khovanskii \cite{Kh77}
for a theorem on generic systems of Laurent polynomials
with fixed Newton polyhedra. Like him, we have to deal 
with many charts, and we use Sard's theorem \cite{Sa42}
that the subset of critical values of a $C^\infty$-map
between $C^\infty$-manifolds has measure 0.
The analogue of the toric compactification which 
Khovanskii constructs is in our situation the space
$\P^\AAA W$.

The product $\P^\AAA W$ had been considered implicitly
in \cite{MM97} and explicitly in \cite{Vi17}, but not
in many other papers on games in normal form,
although the embedding $G\hookrightarrow \P^\AAA W$
is natural.

We do not provide applications of Theorem \ref{t3.1}
in this paper. But we expect that Theorem \ref{t3.1} 
will be useful and that applications will come.

In the case of two-player games, it is not so difficult
to make precise what {\it generic} means and which
properties such games should have \cite{LH64}\cite{St99}.
There one can refer essentially to linear algebra.

For games with three or more players, it is more difficult.
Other proofs, that in a generic game the set 
of Nash equilibria is finite and each Nash equilibrium
is regular, were given in \cite[ch.~5]{Ha73}, 
\cite[Theorem 3]{GPS93} and \cite[Theorem 2]{Ri94}. 
Harsanyi used (as we do) Sard's theorem,
G\"ul, Pearce and Stacchetti used a variant due to
Stacchetti and Reinoza of Sard's theorem,
and Ritzberger used a parametric transversality theorem.

Section \ref{c2} introduces more formally finite
games in normal form and their mixed extensions,
and it sets some notations.
Section \ref{c3} formulates the main result
Theorem \ref{t3.1}.
Section \ref{c4} gives some background material
from differential topology.
Section \ref{c5} proves Theorem \ref{t3.1}.

\section{The mixed extension of a finite game}\label{c2}
\setcounter{equation}{0}

This section introduces finite games, their mixed
extensions, best reply maps and Nash equilibria,
and it sets some notations.

\begin{definition}\label{t2.1}
(a) $(\AAA,S,U)$ denotes a finite game. Here 
$m\in\N=\{1,2,3,...\}$, $\AAA:=\{1,...,m\}$ is the set of
players, $S^i=\{s^i_0,...,s^i_{n_i}\}$ with $n_i\in\N$ 
is the set of pure strategies of player $i\in\AAA$, 
$S=S^1\times ...\times S^m$ is the set of pure strategy
combinations, $U^i:S\to\R$ is the utility function of player
$i$, and $U=(U^1,...,U^m):S\to\R^m$. 
We denote $J^i:=\{1,...,n_i\}$, $J^i_0:=\{0\}\cup J^i$
and $J_0^\AAA:=\prod_{i=1}^mJ^i_0$.
The pure strategy combinations are given as tuples 
$(s^1_{j_1},...,s^m_{j_m})\in S$ with 
$\uuuu{j}=(j_1,...,j_m)\in J_0^\AAA$. 

\medskip
(b) $(\AAA,G,V)$ denotes the {\it mixed extension}
of the finite game in (a). Here 
\begin{eqnarray*}
W^i:=\bigoplus_{j=0}^{n_i}\R\cdot s^i_j,&&
W:=W^1\times ...\times W^m,\\
A^i:=\{\sum_{j=0}^{n_i}\gamma^i_js^i_j\in W^i\,|\, 
\sum_{j=0}^{n_i}\gamma^i_j=1,\},&&
A:=A^1\times ...\times A^m\subset W,\\
G^i:=\{\sum_{j=0}^{n_i}\gamma^i_js^i_j\in A^i\,|\, 
\gamma^i_j\in[0,1]\},&&
G:=G^1\times ...\times G^m\subset A.
\end{eqnarray*} 
So, $W^i$ and $W$ are real vector spaces, 
$A^i\subset W^i$ and $A\subset W$ are
affine linear subspaces of codimension $1$ respectively $m$,
$G^i\subset A^i$ is a simplex in $A^i$ of the same
dimension $n_i$ as $A^i$, and $G\subset A$
is a product of simplices, so especially a convex
polytope, and it has the same 
dimension $\sum_{i=1}^m n_i$ as $A$. 
The map $V^i_W:W\to\R^m$ is the multilinear extension of $U^i$,
\begin{eqnarray}\label{2.1}
V^i_W(g)&:=&\sum_{(j_1,...,j_m)\in J_0^\AAA}
\Bigl(\prod_{k=1}^m\gamma^k_{j_k}\Bigr)
\cdot U^i(s^1_{j_1},...,s^m_{j_m}),\\
\textup{where}\quad g&=&(g^1,...,g^m)\in W\textup{ with }
g^k=\sum_{j=0}^{n_k}\gamma^k_js^k_j,\nonumber
\end{eqnarray}
$V^i_A:A\to\R$ is the restriction of $V_W$ to
$A$, and $V^i:G\to\R$ is the restriction of 
$V^i_W$ to $G$. Then $V=(V^1,...,V^m):G\to\R^m$. 
An element $g\in G$ is called a 
{\it mixed strategy combination}.
The support of an element 
\begin{eqnarray*}
g^i=\sum_{j=0}^{n_i}\gamma^i_js^i_j\in W^i
\quad\textup{is the set} \quad
\supp(g^i):=\{j\in J^i_0\,|\,\gamma^i_j\neq 0\}.
\end{eqnarray*}
We also denote 
$G^{-i}:=G^1\times ...\times G^{i-1}\times G^{i+1}\times ...
\times G^m$, and its elements
$g^{-i}:=(g^1,...,g^{i-1},g^{i+1},...,g^m)\in G^{-i}$. 
We follow the standard (slightly incorrect)
convention and identify $G^i\times G^{-i}$ with $G$
and $(g^i,g^{-i})$ with $g$. 

\medskip
(c) Fix $i\in\AAA$. The {\it best reply map} 
$r^i:G^{-i}\to\PP(G^i)$ associates to each element
$g^{-i}\in G^{-i}$ the set of its best replies in $G^i$,
\begin{eqnarray}\label{2.2}
r^i(g^{-i}):=\{g^i\in G^i\,|\, V^i(\www g^i,g^{-i})
\leq V^i(g^i,g^{-i})\textup{ for any }\www g^i\in G^i\}.
\end{eqnarray}
Its graph is the set 
$$\Gr(r^i):=\bigcup_{g^{-i}\in G^{-i}}r^i(g^{-i})\times 
\{g^{-i}\}\subset G^i\times G^{-i}=G.$$ 
A {\it Nash equilibrium}
is an element of the set $\NN:=\bigcap_{i\in \AAA}\Gr(r^i)$. 
The set $\NN$ is the set of Nash equilibria. 

\medskip
(d) Write $\uuuu{\gamma}^i=(\gamma^i_1,...,\gamma^i_{n_i})$,
\begin{eqnarray*}
\uuuu{\gamma}&:=&(\uuuu{\gamma}^1;...;\uuuu{\gamma}^m) 
=(\gamma^1_1,...,\gamma^1_{n_1};...;\gamma^m_1,...,
\gamma^m_{n_m})\quad\textup{and}\\
\uuuu{\gamma}^{-i}&:=&(\uuuu{\gamma}^1;...;\uuuu{\gamma}^{i-1};
{\uuuu\gamma}^{i+1};...,\uuuu{\gamma}^m).
\end{eqnarray*}
Fix a subset $\BB\subset\AAA$. 
The monomials $\prod_{i\in \BB}\gamma^i_{j_i}$ for 
$(j_i\,|\, i\in \BB)\in \prod_{i\in \BB}J^i_0$ in 
$\R[\gamma^i_0,\uuuu{\gamma}^i\,|\, i\in \BB]$
(with $\prod_{i\in\emptyset}(..):=1$)
are called {\it $\BB$-multilinear}. A polynomial which is a
real linear combination of $\BB$-multilinear monomials is
also called $\BB$-multilinear. Especially $V^i_W$ is
$\AAA$-multilinear. 
A polynomial in $\R[\gamma^i_0,\uuuu{\gamma}^i\,|\, i\in \AAA]$
is {\it multi affine linear} if each monomial in it with 
nonvanishing coefficient is $\CC$-multilinear for a suitable
set $\CC\subset \AAA$. 
\end{definition}

The set $\NN$ of Nash equilibria is not empty. 
This was first proved by Nash \cite{Na51}.
The following lemma is rather trivial, but worth to be noted.

\begin{lemma}\label{t2.2}
Let $(\AAA,G,V)$ be the mixed extension of a finite game
$(\AAA,S,U)$. 

(a) The tuple 
$\uuuu{\gamma}^i$ is a tuple 
of (affine linear) coordinates on $A^i$, because in 
$A^i$ we have $\gamma^i_0=1-\sum_{j=1}^{n_i}\gamma^i_j$.
The tuple 
$\uuuu{\gamma}=(\uuuu{\gamma}^1;...;\uuuu{\gamma}^m)$ is a tuple 
of (affine linear) coordinates on $A$.
The map $V^i_A$ is a multi affine linear
polynomial in $\uuuu{\gamma}$. It has the shape
\begin{eqnarray}\label{2.3}
V^i_A(g)=\kappa^i(\uuuu{\gamma}^{-i})
+\sum_{j=1}^{n_i}\gamma^i_j\cdot
\lambda^i_j(\uuuu{\gamma}^{-i})\quad\textup{for}\quad
g\in A, 
\end{eqnarray}
where $\kappa^i$ and all $\lambda^i_j$ are unique
multi affine linear polynomials in $\uuuu{\gamma}^{-i}$.
Define additionally
\begin{eqnarray}\label{2.4}
\lambda^i_0:=0\quad\textup{for }i\in\AAA.
\end{eqnarray}

(b) An element $g=(g^i,g^{-i})\in G$ is in $\Gr(r^i)$ if
and only if the following holds. 
\begin{eqnarray}
\lambda^i_j(\uuuu{\gamma}^{-i})-\lambda^i_k(\uuuu{\gamma}^{-i})
=0&&\textup{if }j,k\in \supp(g^i),\label{2.5}\\
\lambda^i_j(\uuuu{\gamma}^{-i})-\lambda^i_k(\uuuu{\gamma}^{-i})
\geq 0&&\textup{if }j\in \supp(g^i),k\notin\supp(g^i).\label{2.6}
\end{eqnarray}
\end{lemma}

{\bf Proof:}
(a) Part (a) holds because $A^i$ is the affine hyperplane 
in $W^i$ defined by $\gamma^i_0=1-\sum_{j=1}^{n_i}\gamma^i_j$. 

(b) If $g\in\textup{Gr}(r^i)$ then a change of 
$g^i=\sum_{j=0}^{n_i}\gamma^i_js^i_j$
may not increase $V^i_A(g)$ in \eqref{2.3}. 
Therefore $\lambda^i_j(\uuuu{\gamma}^{-i})$ for 
$j\in\supp(g^i)$ must be the maximum of all 
$\lambda^i_k(\uuuu{\gamma}^{-i})$. 
This includes the case $j=0$ if $0\in\supp(g^i)$, and
it includes the case $k=0$ if $0\notin\supp(g^i)$. 
\hfill$\Box$

\section{Compactification of $A$ and
generic games}\label{c3}
\setcounter{equation}{0}

Consider as in section \ref{c2} a finite set
$\AAA=\{1,...,m\}$ of players and for each player
$i\in\AAA$ a finite set $S^i=\{s^i_0,...,s^i_{n_i}\}$ 
with $n_i\in\N$ of pure strategies. 
Let $\UU^i=\R^S$ be the set of all possible
utility functions $U^i:S\to\R$. The set of all possible 
utility maps $U=(U^1,...,U^m)$ is then 
$\UU:=\prod_{i=1}^m\UU^i\cong (\R^S)^m$. 

Consider a fixed map $U$, the tuple of 
all hyperplanes in $A$ which bound $G\subset A$
and the subvarieties
$(\lambda^i_j-\lambda^i_k)^{-1}(0)\subset A$
for $i\in \AAA$ and $j,k\in J^i_0$ with $j<k$. 
By Lemma \ref{t2.2} (b), 
the graphs $\Gr(r^i)$ of the best
reply maps and the set $\NN$ of Nash equilibria 
are determined by the geometry of these hyperplanes and
these subvarieties.

The hyperplanes which bound $G$ are smooth and transversal.
Theorem \ref{t3.1} below will imply that for generic $U\in\UU$
also the subvarieties $(\lambda^i_j-\lambda^i_k)^{-1}(0)$
are smooth hypersurfaces in $A$ and that they and the
hyperplanes in $A$ which bound $G$ 
are as transversal as possible. 

The fact that $V^i_W$ is $\AAA$-multilinear motivates to
consider the natural compactification of $A^i$ to the
real projective space $\P W^i$, which is the set 
$(W^i-\{0\})/\R^*$ of lines through 0 in $W^i$, 
and the natural compactification of $A$ to
the product of real projective spaces
\begin{eqnarray}\label{3.1}
\P^\AAA W:=\prod_{i=1}^m\P W^i.
\end{eqnarray}
We denote $\P^{-i}W:=\prod_{j\in\AAA-\{i\}}\P W^i$,
and we identify (following the slightly incorrect convention
in Definition \ref{t2.1} (b)) $\P W^i\times \P^{-i}W$ 
with $\P^\AAA W$. Under the projection 
\begin{eqnarray}\label{3.2}
pr_W:\prod_{i=1}^m(W^i-\{0\})\to\P^\AAA W,
\end{eqnarray}
the affine linear space $A\subset W$ embeds
as a Zariski open subset into $\P^\AAA W$. The complement is
the union of $m$ hyperplanes
\begin{eqnarray}\label{3.3}
H^{i,\infty}:=(\P W^i-A^i)\times \P^{-i}W 
\subset \P W^i\times \P^{-i}W=\P^\AAA W
\end{eqnarray}
for $i\in\AAA$. 
For a subset $B^i\subset A^i$, let $\oooo{B^i}^{Zar}$ 
denote its Zariski closure in $\P W^i$,
that is the smallest real algebraic subvariety in $\P W^i$
which contains $B^i$. 
The Zariski closure in $\P^\AAA W$ of a subset 
$B\subset\P^\AAA W$ is denoted by $\oooo{B}^{Zar}$.
For $i\in\AAA$ and $j\in J^i_0$ denote by 
\begin{eqnarray}\label{3.4}
H^{i,j}:=\oooo{\{g^i\in A^i\,|\, \gamma^i_j=0\}}^{Zar}
\times \P^{-i} W\subset \P W^i\times \P^{-i}W=\P^\AAA W
\end{eqnarray}
the Zariski closures in $\P^\AAA W$ of the hyperplanes
in $A$ which bound $G$. 
Define 
$$J^{i,2}_0:=\{(j,k)\in J^i_0\times J^i_0\,|\, j<k\}.$$
For $i\in\AAA$ and $(j,k)\in J^{i,2}_0$ 
consider the difference $\lambda^i_j-\lambda^i_k$
as a function on $A$ (so lift it as a function from
$\prod_{l\neq i}A^l$ to $A$), and consider
the Zariski closure in $\P^\AAA W$
\begin{eqnarray}\label{3.5}
H^{i,(j,k)}:=\oooo{(\lambda^i_j-\lambda^i_k)^{-1}(0)}^{Zar}
\subset\P^\AAA W. 
\end{eqnarray}

The notion {\it everywhere transversal} in Theorem \ref{t3.1} 
is defined in Definition \ref{t4.4} (b). 
In Theorem \ref{t3.1}, a subset
of the set of all hyperplanes in \eqref{3.3} and \eqref{3.4} 
and all subvarieties in \eqref{3.5} is considered. 
Such a subset is characterized by its set of indices,
namely sets $T^i\subset J^i_0\cup\{\infty\}$ and sets 
$R^i\subset J^{i,2}_0$ for $i\in\AAA$ define a subset 
\begin{eqnarray}\label{3.6}
\bigcup_{i\in\AAA}\Bigl(\{H^{i,j}\,|\, j\in T^i\}\cup 
\{H^{i,(j,k)}\,|\, (j,k)\in R^i\}\Bigr).
\end{eqnarray}
A set $R^i\subset J^{i,2}_0$ defines
a graph with vertex set $J^i_0$ and set $R^i$ of edges. 
The set in \eqref{3.6} is called {\it good} if the graph
$(J^i_0,R^i)$ is a union of trees. 

The reason for the introduction of this notion of 
a {\it good} set is that in the case of $|\supp(g^i)|\geq 3$
there is some redundancy in the equations \eqref{2.5}.
Equality for all pairs $(j,k)$ with $j,k\in \supp(g^i)$ is implied
by equality for a set of pairs $(j,k)$ such that the graph with
vertex set $\supp(g^i)$ and edge set this set of pairs is a tree.

\begin{theorem}\label{t3.1}
Let $\AAA=\{1,...,m\}$ and $S$ be as in section \ref{c2}
and as above. 
There is a semialgebraic subset $\DD\subset\UU$ of codimension
at least one (equivalently, it is of Lebesgue measure 0)
such that for any tuple $U\in\UU-\DD$ 
of utility functions the following holds.
The hyperplanes $H^{i,j},i\in\AAA,j\in J^i_0\cup\{\infty\}$,
and the subvarieties 
$H^{i,(j,k)},i\in\AAA,(j,k)\in N^ {i,2}_0$, 
are smooth hypersurfaces in $\P^\AAA W$.
Any good subset of them is everywhere transversal.
\end{theorem}

The proof of Theorem \ref{t3.1} will be given in section 
\ref{c5}. Before, section \ref{c4} will recall some basic
notions and facts from differential topology.

\section{Transversality of submanifolds}\label{c4}
\setcounter{equation}{0}

This section recalls two basic facts from differential
topology, Sards theorem and the implicit function theorem.
It defines transversality of submanifolds and it formulates
a useful lemma.

\begin{definition}\label{t4.1}
Let $M$ and $N$ be $C^\infty$-manifolds, and let 
$F:M\to N$ be a $C^\infty$-map. A point $p\in M$ is a
{\it regular point} of $F$ if the linear map
$dF|_{T_pM}:T_pM\to T_{F(p)}N$ is surjective.
A point $p\in M$ is a {\it critical point} if it is not
a regular point. A point $q\in N$ is a {\it regular value}
if either $F^{-1}(q)=\emptyset$ or every point
$p\in F^{-1}(q)$ is a regular point. A point $q\in N$
is a {\it critical value} if it is not a regular value.
\end{definition}

The following theorem is famous.
It is also crucial in the proof of Theorem \ref{t3.1}.

\begin{theorem}\label{t4.2} (Sard's theorem \cite{Sa42},
see also \cite[Theorem 2.1]{BL75})
Let $F:M\to N$ be a $C^\infty$-map between $C^\infty$-manifolds.
The subset of $N$ of critical values of $F$ has Lebesgue
measure 0.
\end{theorem}

In our situation, the setting is semialgebraic. Therefore
the subset of critical values is then a semialgebraic
subset of $N$ of Lebesgue measure 0. This means that it has
everywhere smaller dimension than $N$.  
A variant of Sard's theorem in the semialgebraic setting
says precisely this \cite[ch. 9.5]{BCR98}\cite[2.5.12]{BR90}.

The implicit function theorem says how a map 
$F:M\to N$ looks near a regular point $p\in M$.
It is a trivial fibration with smooth fibers.

\begin{theorem}\label{t4.3} 
(Implicit function theorem, e.g. \cite[Theorem 1.3]{BL75})
Let $F:M\to N$ be a $C^\infty$-map between 
$C^\infty$-manifolds, and let $p\in M$ be a regular point
of $F$. Then $\dim M\geq \dim N$, and there are open neighborhoods $U_1\subset M$ of $p$
and $U_2\subset N$ of $F(p)$ with $U_2\supset F(U_1)$, 
open balls $B_1\subset\R^{\dim M}$
around $0$ and $B_2\subset \R^{\dim N}$ around $0$ and
$C^\infty$-diffeomorphisms $\varphi_1:B_1\to U_1$ and
$\varphi_2:B_2\to U_2$ with the following property.
The map $\varphi_2^{-1}\circ F\circ\varphi_1:B_1\to B_2$ 
is the standard projection in \eqref{4.1},
\begin{eqnarray}\label{4.1}
\begin{array}{ccc}
U_1 & \stackrel{F}{\longrightarrow} & U_2\\
\hspace*{0.5cm}\uparrow\varphi_1 & & 
\hspace*{0.5cm}\uparrow\varphi_2\\
B_1 & \stackrel{\varphi_2^{-1}\circ F\circ\varphi_1}{\longrightarrow} & B_2\\
(x_1,...,x_{\dim M})&\mapsto& (x_1,...,x_{\dim N})
\end{array}
\end{eqnarray}
\end{theorem}

\begin{definition}\label{t4.4}
(a) Let $M$ be a $C^\infty$-manifold, let $p\in M$ and
let $H\subset M$ be a subset with $p\in H$ and which is in a 
neighborhood of $p$ a $C^\infty$-submanifold of $M$.
A {\it defining map} $F$ for the pair $(H,p)$ is
a function $F:U\to\R^n$ with $U\subset M$ an open neigborhood 
of $p$ such that $F$ is regular on each point of $U$ and
$H\cap U=F^{-1}(F(p))$. Then $n=\dim M-\dim H\cap U$. 

Remark: Defining maps for $(H,p)$ exist and are related
by local diffeomorphisms.

(b) Let $M$ be a $C^\infty$-manifold, and let
$H_1,...,H_a$ be $C^\infty$-submanifolds.
Now it will be defined when they are 
{\it transversal at a point $p\in M$}.
\begin{list}{}{}
\item[(i)] They are transversal at $p\in M-\cup_{i=1}^aH_i$.
\item[(ii)] They are transversal at $p\in \cap_{i=1}^aH_i$
if for some (or for any, that is equivalent) tuple of 
defining maps $F_i:U\to \R^{n_i}$ of $(H_i,p)$
($i\in\{1,...,a\}$) with joint definition domain $U$ 
the map $(F_1,...,F_a):U\to\R^{n_1+...+n_a}$ is regular at $p$.
\item[(iii)] They are transversal at $p\in \cup_{i=1}^aH_i$
if the subset $\{H_i\,|\, p\in H_i\}$ of the set of manifolds 
$H_1,...,H_a$ is transversal at $p$
(this was defined in part (ii)).
\end{list}
Finally, they are {\it transversal} (or {\it everywhere
transversal}) if they are transversal at each point of $M$.
\end{definition}

Part (a) of the following lemma states an obvious consequence
on the intersection of a family of transversal submanifolds.
Part (b) gives a useful criterion for proving
transversality of a family of submanifolds.

\begin{lemma}\label{t4.5}
Let $M$ be a $C^\infty$-manifold, and let $H_1,...,H_a$ 
be $C^\infty$-submanifolds. 

(a) If they are transversal then either 
$\cap_{i=1}^aH_i=\emptyset$ or this intersection 
is a submanifold of $M$ of dimension
$\dim M-\sum_{i=1}^a\codim H_i$. Especially, in the second
case, this number is non-negative. 

(b) Suppose that for some $b\in\{1,..,a-1\}$ the submanifolds
$H_1,...,H_b$ are transversal at a point $p\in \cap_{i=1}^aH_i$.
The intersection $L:=\cap_{i=1}^bH_i$ is by part (a) 
in a suitable open neighborhood of $p$ a submanifold of $M$. 
The following two conditions are equivalent:
\begin{list}{}{}
\item[(i)]
$H_1,...,H_a$ are transversal at $p$.
\item[(ii)]
For $j\in\{b+1,...,a\}$ there are defining maps
$F_j:U_j\to \R^{n_j}$ on open neighborhoods $U_j\subset M$ of $p$
such that for $U:=\bigcap_{j=b+1}^aU_j$ the intersection
$U\cap L$ is a submanifold of $U$ and the map
$(F_{b+1},...,F_a)|_{U\cap L}:U\cap L\to\R^{n_{b+1}+...+n_a}$
is regular at $p$. 
\end{list}
\end{lemma}

{\bf Proof:}
(a) Part (a) follows immediately from the definition
of transversality and the implicit function theorem.

(b) By the implicit function theorem we can choose 
an open neighborhood $U$ of $p$ in $M$, coordinates
$\uuuu{x}=(x_1,...,x_m)$ on $U$  
with $\uuuu{x}(p)=0$ and defining maps 
$F_i:=U\to\R^{n_i}$ of $(H_i,p)$ for $i\in\{1,...,a\}$
such that $(F_1,...,F_b)=(x_1,...,x_{\sum_{i=1}^bn_i})$.
Condition (i) is equivalent to 
\begin{eqnarray}\label{4.2}
\rank\begin{pmatrix}\frac{\paa F_i}{\paa x_j}(p)
\end{pmatrix}_{\tiny\begin{array}{ll}i\in\{1,...,a\}\\j\in\{1,...,\dim M\}\end{array}}=\sum_{i=1}^a n_i.
\end{eqnarray}
Condition (ii) is equivalent to
\begin{eqnarray}\label{4.3}
\rank\begin{pmatrix}\frac{\paa F_i|_{U\cap L}}{\paa x_j}(p)
\end{pmatrix}_{\tiny
\begin{array}{ll}i\in\{b+1,...,a\}\\j\in\{1+\sum_{i=1}^bn_i,...,\dim M\}\end{array}}=\sum_{i=b+1}^a n_i.
\end{eqnarray}
The matrix in \eqref{4.2} has a block triangular shape
$\left(\begin{array}{c|c} {\bf 1} & 0 \\ \hline * & * 
\end{array}\right)$
with the matrix in \eqref{4.3} in the lower right place.
Therefore \eqref{4.2} and \eqref{4.3} are equivalent. 
\hfill$\Box$

\begin{remarks}\label{t4.6}
(i) In condition (ii) in Lemma \ref{t4.5} (b), one can 
replace the existence of the defining maps by demanding
that the condition holds for all choices of defining maps. 

(ii) In the proof in section \ref{c5} of the transversality 
in Theorem \ref{t3.1}, Lemma \ref{t4.5} (b) will be useful. 
The hyperplanes $H^{i,j}$ in a given set in \eqref{3.6} 
take the role of the submanifolds $H_1,...,H_b$. 
The reason is that they are fixed if one moves $U\in\UU$,
while the other subvarieties 
$H^{i,(j,k)}$ move if one moves $U\in\UU$.

(iii) In Lemma \ref{t4.5} (b) one cannot replace condition (ii)
by the simpler condition (ii)': 
{\it $L\cap H_{b+1},...,L\cap H_a$ are submanifolds and are
transversal at $p$.} The reason is that $L\cap H_j$ for
some $j\in\{b+1,...,a\}$ can be a submanifold of the
correct dimension although $L$ and $H_j$ are not transversal. 
\end{remarks}

\section{Proof of Theorem \ref{t3.1}}\label{c5}
\setcounter{equation}{0}

This section is devoted to the proof of Theorem \ref{t3.1}.

First we discuss affine charts on $\P W^i$ and on 
$\P^\AAA W$.
The linear coordinates $(\gamma^i_0,...,\gamma^i_{n_i})$ on
$W^i$ induce $n_i+1$ affine charts of the projective
space $\P W^i$. But none of them contains the whole set $G^i$.
The following linear coordinates 
$\uuuu{\www\gamma}^i:=(\www\gamma^i_0,...,\www\gamma^i_{n_i})$ 
on $W^i$ are equally natural, and one of the affine charts
which they induce on $\P W^i$ will turn out to be $A^i$,
\begin{eqnarray}\label{5.1}
\www\gamma^i_j:=\gamma^i_j\textup{ for }j\in J^i,\quad
\www\gamma^i_0:=\sum_{j=0}^{n_i}\gamma^i_j.
\end{eqnarray}
From now on we use the induced homogeneous coordinates
$(\www\gamma^{i_0}:...:\www\gamma^i_{n_i})$ on $\P W^i$
if not said otherwise.
The $n_i+1$ induced affine charts of $\P W^i$ are 
\begin{eqnarray}\label{5.2}
A^i_j:=\{(\www\gamma^i_0:...:\www\gamma^i_{n_i})
\in\P W^i\,|\, \uuuu{\www\gamma}^i\in W^i,\www\gamma^i_j=1\}
\subset\P W^i
\quad\textup{for }j\in J^i_0.
\end{eqnarray}
$A^i_j$ comes equipped with natural coordinates
$\uuuu{\gamma}^{i,j}=(\gamma^{i,j}_0,...,\gamma^{i,j}_{j-1},
\gamma^{i,j}_{j+1},...,\gamma^{i,j}_{n_i})$ from the 
isomorphism
\begin{eqnarray*}
\R^{n_i}\to A^i_j,\quad\uuuu{\gamma}^{i,j}\mapsto
(\gamma^{i,j}_0:...:\gamma^{i,j}_{j-1}:1:
\gamma^{i,j}_{j+1}:...:\gamma^{i,j}_{n_i}),
\end{eqnarray*}
and we have a natural embedding
\begin{eqnarray*}
\alpha^{i,j}:A^i_j\hookrightarrow W^i,\quad 
\uuuu{\gamma}^{i,j}\mapsto 
(\gamma^{i,j}_0,...,\gamma^{i,j}_{j-1},1,
\gamma^{i,j}_{j+1},...,\gamma^{i,j}_{n_i}),
\end{eqnarray*}
with $\pr^i_W\circ\alpha^{i,j}=\id$ where
$\pr^i_W:W^i-\{0\}\to \P W^i$ is the natural projection.
The image $\alpha^{i,0}(A^i_0)\subset W^i$ 
coincides with the set $A^i\subset W^i$ 
in Definition \ref{t2.1} (b), and the identification
$\uuuu{\gamma}^{i,0}=\uuuu{\gamma}^i$ identifies $A^i_0$
with $A^i$. 

All possible products over $i\in\AAA$ of these charts
give the affine charts
\begin{eqnarray}\label{5.3}
A_{\uuuu{j}}:=\prod_{i=1}^mA^i_{j_i}\cong\R^{\sum_{i=1}^mn_i}
\quad\textup{for }
\uuuu{j}=(j_1,...,,j_m)\in J_0^\AAA=\prod_{i=1}^mJ^i_0.
\end{eqnarray}
of $\P^\AAA W$ with coordinates $\uuuu{\gamma}^{\uuuu{j}}
=(\uuuu{\gamma}^{1,j_1};...;\uuuu{\gamma}^{m,j_m})$. 
The embeddings 
$\alpha^{i,j_i}:A^i_{j_i}\hookrightarrow W^i$
combine to an embedding 
\begin{eqnarray*}
\alpha^{\uuuu{j}}:A_{\uuuu{j}}\hookrightarrow W
\end{eqnarray*}
with $\alpha^{\uuuu{j}}\circ \pr_W=\id$ on $A_{\uuuu{j}}$. 
The image $\alpha^{\uuuu{0}}(A_{\uuuu{0}})\subset W$ 
coincides with the set $A\subset W$ 
in Definition \ref{t2.1} (b), and the identification
$\uuuu{\gamma}^{\uuuu{0}}=\uuuu{\gamma}$ identifies 
$A_{\uuuu{0}}$ with $A$.

Now we discuss the hyperplanes $H^{i,j}$ and the
subvarieties $H^{i,(j,k)}$ from section \ref{c3}
with respect to the new linear coordinates
$\uuuu{\www\gamma}=(\uuuu{\www\gamma}^1;...;\uuuu{\www\gamma}^m)$
on $W$. 
We define for $i\in\AAA$ 
$$\www{H}^{i,j}:=H^{i,j}\textup{ for  }j\in J^i,\quad
\www{H}^{i,0}:=H^{i,\infty},\quad \www{H}^{i,\infty}:=H^{i,0}.$$
We have for $i\in\AAA$ 
\begin{eqnarray}\label{5.4}
\www{H}^{i,j}=H^{i,j}=\pr_W\Bigl(\{\uuuu{\www\gamma}\in 
\prod_{k=1}^m(W^k-\{0\})\,|\, \www\gamma^i_j=0\}\Bigr)
\quad\textup{for }j\in J^i,\\
\www{H}^{i,0}=H^{i,\infty}=\pr_W\Bigl(\{\uuuu{\www\gamma}\in 
\prod_{k=1}^m(W^k-\{0\})\,|\, \www\gamma^i_0=0\}\Bigr),
\label{5.5}\\
\www{H}^{i,\infty}=H^{i,0}=\pr_W\Bigl(\{\uuuu{\www\gamma}\in 
\prod_{k=1}^m(W^k-\{0\})\,|\, 
\www\gamma^i_0-\sum_{j=1}^{n_i}\www\gamma^i_j=0\}\Bigr).
\label{5.6}
\end{eqnarray}
Obviously for each chart $A_{\uuuu{j}}$, the complement
is
\begin{eqnarray}\label{5.7}
\P^\AAA W-A_{\uuuu{j}}=\bigcup_{i=1}^m \www{H}^{i,j_i}.
\end{eqnarray}
For each of the hyperplanes  
$\www{H}^{i,j}\cap A_{\uuuu{j}}$ with $j\in J^i_0-\{j_i\}$, 
in the chart $A_{\uuuu{j}}$ a defining map 
(in the sense of Definition \ref{t4.4} (a))
is $\www{\gamma}^i_j\circ\alpha^{\uuuu{j}}$. 
For each of the hyperplanes
$\www{H}^{i,\infty}\cap A_{\uuuu{j}}$ in the chart 
$A_{\uuuu{j}}$ a defining map
is $\bigl(\www{\gamma}^i_0-\sum_{j=1}^{n_i}\www{\gamma}^i_j\bigr)
\circ\alpha^{\uuuu{j}}$. 

Obviously, each subset of the set
$\{\www H^{i,j}\,|\, i\in\AAA,j\in J^i_0\cup\{\infty\}\}$
of all these hyperplanes is transversal 
everywhere in $\P^\AAA W$. 

The map $V^i_W:W\to\R$ is an $\AAA$-multilinear map
also in the new linear coordinates 
$\uuuu{\www\gamma}$. Write it as a sum 
\begin{eqnarray}\label{5.8}
V^i_W(\uuuu{\www\gamma})=\www\gamma^i_0\cdot K^i
(\uuuu{\www\gamma}^{-i})+
\sum_{j=1}^{n_i}\www\gamma^i_j\cdot\Lambda^i_j
(\uuuu{\www\gamma}^{-i})
\quad\textup{for}\quad \uuuu{\www\gamma}=
(\uuuu{\www\gamma}^1,...,\uuuu{\www\gamma}^m)\in W.
\end{eqnarray}
Here $K^i(\uuuu{\www\gamma}^{-i})$ and
$\Lambda^i_j(\uuuu{\www\gamma}^{-i})$ are 
$\AAA-\{i\}$-multilinear in $\uuuu{\www\gamma}^{-i}$. 
Define additionally
\begin{eqnarray}\label{5.9}
\Lambda^i_0:=0.
\end{eqnarray}
We have for $i\in\AAA$ and $(j,k)\in J^{i,2}_0$
\begin{eqnarray}\label{5.10}
H^{i,(j,k)}=\pr_W\Bigl(\{\uuuu{\www\gamma}\in 
\prod_{k=1}^m(W^k-\{0\})\,|\, (\Lambda^i_j-\Lambda^i_k)
(\uuuu{\www\gamma})=0\}\Bigr).
\end{eqnarray}
The intersection $H^{i,(j,k)}\cap A_{\uuuu{l}}$ of the 
subvariety $H^{i,(j,k)}\subset\P^\AAA W$ with the affine chart
$A_{\uuuu{l}}$ for $\uuuu{l}=(l_1,...,l_m)\in J_0^\AAA$ is the zero set 
of the function 
$(\Lambda^i_j-\Lambda^i_k)\circ\alpha^{\uuuu{l}}$.

Especially, the identifications 
$\uuuu\gamma^{\uuuu{0}}=\uuuu\gamma$ and 
$A_{\uuuu{0}}=A$ yield the identification
$\Lambda^i_j\circ\alpha^{\uuuu{0}}=\lambda^i_j$. 
This connects the description of $H^{i,(j,k)}$ in \eqref{5.10}
with the one in \eqref{3.5}. 

Now the main point will be to prove that there is a subset 
$\DD\subset\UU$ of Lebesgue measure 0 with the properties 
in Theorem \ref{t3.1}. Later we will argue that it is
semialgebraic. Then $\DD$ having Lebesgue measure 0 is 
equivalent to $\DD$ having codimension at least one in $\UU$.

Choose some $\uuuu{l}\in J_0^\AAA$ 
and consider the chart $A_{\uuuu{l}}$ of $\P^\AAA W$. 
Choose for each $i\in\AAA$ a set
$T^i\subset J^i_0\cup\{\infty\}$ with $l_i\notin T^i$ 
(we exclude the case $l_i$ because 
$\www{H}^{i,l_i}\subset \P^\AAA W-A_{\uuuu{l}}$ by
\eqref{5.7}) 
and a set $R^i\subset J^{i,2}_0$ such that 
$(J^i_0,R^i)$ is a union of trees
and such that (for all $i$ together) 
$\bigcup_{i=1}^m R^i\neq \emptyset$. 
Though $T^i=\emptyset$ (for some or all $i\in\AAA$) and
$R^i=\emptyset$ (for some, but not all $i\in\AAA$) are allowed.
Now we go into Lemma \ref{t4.5} (b) with 
$$M=A_{\uuuu{l}}\quad\textup{and}\quad
\{H_1,...,H_b\}=\{\www{H}^{i,j}\cap A_{\uuuu{l}}\,|\, i\in\AAA,
j\in T^i\}$$ 
(which is empty if all $T^i=\emptyset$) and 
\begin{eqnarray}\label{5.11}
L=\bigcap_{i\in\AAA,j\in T^i}\www H^{i,j}\cap A_{\uuuu{l}},
\end{eqnarray}
which is an affine linear subspace of $A_{\uuuu{l}}$ of 
codimension $\sum_{i=1}^m|T^i|$.
We can choose a part of the coordinates 
$\uuuu{\gamma}^{\uuuu{l}}$ on $A_{\uuuu{l}}$ 
which forms affine linear coordinates on $L$. 
We choose such a part and call it $\uuuu{\gamma}^{\uuuu{l},L}$.

For a moment, fix one map $U\in \UU$ and consider the 
map 
\begin{eqnarray}\label{5.12}
F^{(U)}:L&\to& N\quad\textup{with }N:=\R^{\sum_{i\in\AAA}|R^i|}\\
\uuuu{\gamma}^{\uuuu{l}}&\mapsto&
\bigl(((\Lambda^i_j-\Lambda^i_k)\circ\alpha^{\uuuu{l}})
(\uuuu{\gamma}^{\uuuu{l}})\,|\, i\in\AAA,(j,k)\in R^i\bigr).
\nonumber
\end{eqnarray}
Each entry 
$(\Lambda^i_j-\Lambda^i_k)\circ\alpha^{\uuuu{l}})|_L$
of this map $F^{(U)}$ is a multi affine linear 
function in those coordinates of $\uuuu{\gamma}^{\uuuu{l},L}$ 
which are not in $\uuuu{\gamma}^{i,\uuuu{l}}$. 

By Sard's theorem (Theorem \ref{t4.2}) the set of critical
values of $F^{(U)}$ has Lebesgue measure 0 in $N$.

Now we consider simultaneously all $U\in \UU$. We claim that
the set of $U\in\UU$ such that $0\in N$ is a critical
value of $F^{(U)}$ has Lebesgue measure 0 in $\UU$. 
The argument for this will connect some natural coordinate
system on $\UU$ with the coefficients of the monomials in the
entries of $F^{(U)}$. 

The space of possible maps
$F^{(U)}$ can be identified with the product
\begin{eqnarray}\label{5.13}
\UU^{(I)}:=\prod_{i\in\AAA,(j,k)\in R^i}
\R[\uuuu{\gamma}^{\uuuu{l},L}
\backslash \uuuu{\gamma}^{i,\uuuu{l}}]_{\bf mal}
\end{eqnarray}
where $\R[\uuuu{\gamma}^{\uuuu{l},L}
\backslash \uuuu{\gamma}^{i,\uuuu{l}}]_{\bf mal}$ 
denotes the space of multi affine linear polynomials in those
coordinates in $\uuuu{\gamma}^{\uuuu{l},L}$ 
which are not in $\uuuu{\gamma}^{i,\uuuu{l}}$ 
(here {\bf mal} stands for 
{\it {\bf m}ulti {\bf a}ffine {\bf l}inear}). 
The space of possible constant summands of the entries of 
$F^{(U)}$ can be identified with the space 
$N=\R^{\sum_{i\in\AAA}|R^i|}$. 

There are natural linear maps $p_{(I)}$ and 
$p_{(II)}$ defined by 
\begin{eqnarray}\label{5.14}
p_{(I)}:\UU\to \UU^{(I)},&& p_{(I)}(U)=F^{(U)},\\
p_{(II)}:\UU^{(I)}\to N,&& p_{(II)}(F)=F(0).
\label{5.15}
\end{eqnarray}
Because of the hypothesis that $(J^i_0,R^i)$ is for any
$i\in\AAA$ a union of trees, the maps $p_{(I)}$ and 
$p_{(II)}$ are surjective. 
To see this, recall that 
the $\AAA$-multilinear maps $V^i_W$ are sums of all monomials 
in the variables $\www\gamma^i_j$ with coefficients, which
form a natural linear coordinate system of the $\R$-vector
space $\UU$. For example, the coefficient of the monomial  
$\www\gamma^i_j\prod_{a\in\AAA-\{i\}}\www\gamma^a_{l_a}$ 
in $V^i_W$ is the coefficient of the constant term in the 
multi affine linear function $\Lambda^i_j\circ\alpha^{\uuuu{l}}$.
The entries of $F^{(U)}$ for a fixed $U$ are the 
multi affine linear functions 
$((\Lambda^i_j-\Lambda^i_k)\circ\alpha^{\uuuu{l}})|_L$ for
$i\in\AAA,(j,k)\in R^i$.
The condition that $(J^i_0,R^i)$ is a union of trees takes care
that sufficiently many of the coefficients of the monomials
in the maps $V^i_W$ 
turn up in the maps $p_{(I)}$ and $p_{(II)}$, so that
these maps are surjective. 

For $F\in\ker(p_{(II)})\subset\UU^{(I)}$ denote by 
$C(F)\subset N$ the set of critical values
of $F:L\to N$. It has Lebesgue measure 0 in $N$ by Sard's theorem
(Theorem \ref{t4.2}). 
For each map $F-c:L\to N$ with $c\in N$ the value $0\in N$
is a critical value only if $c\in C(F)$. 
Therefore $0$ is a critical value of $F^{(U)}:L\to N$ 
for $U\in\UU$ only if
\begin{eqnarray}\label{5.16}
U\in\bigcup_{F\in \ker(p_{(II)})}\bigcup_{c\in C(F)}
p_{(I)}^{-1}(F-c),
\end{eqnarray}
which is a set of Lebesgue measure 0 in $\UU$. 
(see e.g. Theorem 2.8 (Fubini) in \cite{BL75}). 

There are only finitely many charts $A_{\uuuu{l}}$ and
finitely many choices of sets $T^i$ and $R^i$ as above. 
Each leads only to a set of Lebesgue measure 0 in $\UU$.
Also their union $\DD$ has Lebesgue measure 0 in $\UU$. 
If  $U\in\UU-\DD$, then for each map $F^{(U)}:L\to N$
as above, the value $0\in N$ is a regular value. 

For an arbitrary chart $A_{\uuuu{l}}$, 
the special choice $T^a=\emptyset$ for all $a\in \AAA$, 
$R^i=\{(j,k)\}$ for one $i\in \AAA$, and 
$R^a=\emptyset$ for all $a\neq i$ shows that 
$H^{i,(j,k)}\cap A_{\uuuu{l}}$
is a smooth hypersurface for $U\in\UU-\DD$.
It also shows that 
the map $(\Lambda^i_j-\Lambda^i_k)\circ\alpha_{\uuuu{l}}$
is a defining map for the pair $(H^{i,(j,k)},p)$ 
at each point $p\in H^{i,(j,k)}\cap A_{\uuuu{l}}$
(in the sense of Definition \ref{t4.4} (a)). 

For an arbitrary chart $A_{\uuuu{l}}$ and all choices
of the sets $T^i$ and $R^i$, the construction of $\DD$ above
together with Lemma \ref{t4.5} (b) and
Definition \ref{t4.4} (b) show that 
for $U\in\UU-\DD$ any good subset of the smooth hypersurfaces
$H^{i,j}$ and $H^{i,(j,k)}$ is transversal on $A_{\uuuu{l}}$.
Considering all charts together, we obtain all
statements of Theorem \ref{t3.1} except that $\DD$ is 
semialgebraic. 

$\DD$ is semialgebraic because of the following.
For a chart $A_{\uuuu{l}}$ and sets $T^i$ and $R^i$ as above, 
the maps $F^{(U)}:L\to N$ unite to an algebraic map 
\begin{eqnarray}\label{5.17} 
\FF:\UU\times L\to \UU\times N,\quad(U,x)\mapsto (U,F^{(U)}(x)).
\end{eqnarray}
The set of its critical points in $\UU\times L$ with 
critical value $(U,0)\in\UU\times\{0\}\subset\UU\times N$ 
is an algebraic subset of $\UU\times L$. 
Its image under the projection to $\UU$ is semialgebraic.
Also the union $\DD$ of these sets over all choices of charts
$A_{\uuuu{l}}$ and sets $T^i$ and $R^i$ is semialgebraic.
\hfill$\Box$

\begin{remark}\label{t5.1}
This proof of Theorem \ref{t3.1} is inspired by the 
proof of Khovanskii of a theorem on generic systems
of Laurent polynomials with fixed Newton polyhedra. 
He considers complex coefficients, we consider real 
coefficients. But apart from that, our situation and proof
can be seen as a special case of his situation and proof.
It is the main theorem in $\S 2$ in \cite{Kh77}.
Though, he uses, but does not formulate the Lemma \ref{t4.5}.
The analogue of the toric compactification which he has
to construct is in our situation the space $\P^\AAA W$. 
\end{remark}

The proof of Theorem \ref{t3.1} in this section \ref{c5}
implies the following corollary. The regularity which
it expresses is most useful in the case of points 
(e.g. Nash equilibria) in the standard affine chart 
$A_{\uuuu{0}}$ and written with the
standard defining maps in this chart.

\begin{corollary}\label{t5.2}
Consider the situation in Theorem \ref{t3.1}.
Let $U$ be the utility map of a game with $U\in\UU-\DD$
(so it is generic).
Let $A_{\uuuu{l}}\subset\P^\AAA W$ be any one of the
affine charts of $\P^\AAA W$.
Let $\uuuu{\gamma}^{\uuuu{l}}$ be any point in 
$A_{\uuuu{l}}$. Consider any set of smooth hypersurfaces
as in \eqref{3.6} which is good and which contains
the point $\uuuu{\gamma}^{\uuuu{l}}$.
Then the Jacobian of the defining maps in this chart
(which were used in the proof above) for these
hypersurfaces is nondegenerate.
\end{corollary}

{\bf Proof:}
By Theorem \ref{t3.1} the subvarieties in a set as
in \eqref{3.6} are smooth hypersurfaces and everywhere
transversal. The transversality in any affine chart
$A_{\uuuu{l}}$ was shown by proving that the point
$\uuuu{\gamma}^{\uuuu{l}}$ is a regular point of
the tuple of defining maps in this chart 
of the hypersurfaces.
\hfill$\Box$

\end{document}